\documentclass[11pt]{amsart}
\usepackage{graphicx}
\usepackage{amsmath, amsfonts, amssymb, amsthm}

\usepackage{fullpage, color}

\newtheorem{theorem}{Theorem}[section]
\newtheorem{conjecture}[theorem]{Conjecture}

\newtheorem{lemma}[theorem]{Lemma}

\allowdisplaybreaks[4] 
\begin{document}


\title{A proof of some conjectures of Mao on partition rank inequalities}

\author{Ethan Alwaise}
\address{Emory University}
\email{ealwais@emory.edu}

\author{Elena Iannuzzi}
\address{Vassar College}
\email{eliannuzzi@vassar.edu}

\author{Holly Swisher}
\address{Department of Mathematics, Oregon State University, 368 Kidder Hall, Corvallis, OR 97331}
\email{swisherh@math.oregonstate.edu}

\thanks{This work was supported by the National Science Foundation REU Grant DMS-1359173.}

\date{\today}



\begin{abstract}
Based on work of Atkin and Swinnerton-Dyer on partition rank difference functions, and more recent work of Lovejoy and Osburn, Mao has proved several inequalities between partition ranks modulo $10$, and additional results modulo $6$ and $10$ for the $M_2$ rank of partitions without repeated odd parts.  Mao conjectured some additional inequalities.  We prove some of Mao's rank inequality conjectures for both the rank and the $M_2$ rank modulo $10$ using elementary methods.  
\end{abstract}

\maketitle

\section{Introduction and Statement of Results}\label{intro}

For a positive integer $n$, a {\it partition} of $n$ is a non-increasing sequence of positive integers that sum to $n$, where each summand is called a {\it part}. The partition function $p(n)$ counts the number of partitions of $n$, and we define $p(0)=1$.  

The celebrated Ramanujan congruences demonstrate compelling divisibility properties for $p(n)$,
\begin{align*}
p(5n+4) & \equiv 0 \pmod 5 \\
p(7n+5) & \equiv 0 \pmod 7  \\
p(11n+6) & \equiv 0 \pmod{11}.
\end{align*}

Dyson \cite{Dyson} defined the {\it rank} of a partition $\lambda$ to be $l(\lambda)-n(\lambda)$, where $l(\lambda)$ and $n(\lambda)$ denote the largest part and number of parts of $\lambda$, respectively.  Dyson conjectured that this gave a combinatorial explanation for the Ramanujan congruences modulo $5$ and $7$.   In particular, if $N(s,m,n)$ is defined to be the number of partitions of $n$ that have rank congruent to $s$ modulo $m$, then Dyson conjectured that for each residue class $s$,
\begin{align}
N(s,5,5n+4) &= \frac{p(5n+4)}{5}  \label{Dyson1}\\
N(s,7,7n+5) &= \frac{p(7n+5)}{7} \label{Dyson2}.
\end{align}

Atkin and Swinnerton-Dyer \cite{ASD} proved (\ref{Dyson1}), (\ref{Dyson2}) by obtaining generating functions for rank differences of the form $N(s,\ell,\ell n+b)-N(t,\ell,\ell n+b)$ for $\ell=5,7$, and showing that the relevant differences were always $0$ in the setting $(\ell,b)\in\{(5,4),(7,5)\}$.  They determined all of the generating functions for $N(s,\ell,\ell n+b)-N(t,\ell,\ell n+b)$ where $\ell=5,7$, and obtained several interesting identities for the non-Ramanujan cases.  

Lovejoy and Osburn \cite{LO1, LO3, LO2} used similar techniques to obtain interesting generating function representations for rank differences of overpartitions, as well as partitions without repeated odd parts.  For example, let $\lambda$ be a partition without repeated odd parts.  The {\it $M_2$ rank} of $\lambda$ is defined to be
\[
\left\lceil{\frac{l(\lambda)}{2}}\right\rceil - n(\lambda).
\]
Let $N_2(s,m,n)$ count the number of partitions of $n$ with distinct odd parts and $M_2$ rank congruent to $s$ modulo $m$.  Lovejoy and Osburn \cite{LO3} obtained generating function identities for rank differences of the form $N_2(s,\ell,\ell n + b) - N_2(t,\ell,\ell n + b)$ for $\ell = 3$ and $\ell = 5$.  

Most recently, Mao \cite{Mao1, Mao2} has derived generating function formulas for Dyson's rank on partitions modulo $10$, and the $M_2$ rank on partitions without repeated odd parts modulo $6$ and $10$.    In this work he proves a number of inequalities, including for example
\begin{align*}
N(0,10,5n+1)  & > N(4,10,5n+1), \\
N_2(0,6,3n) + N_2(1,6,3n) & > N_2(2,6,3n) + N_2(3,6,3n).
\end{align*} 
Mao gives the following conjectures based on computational evidence.  The first, is for Dyson's rank on unrestricted partitions. 
\begin{conjecture}
Computational evidence suggests that
\begin{align}
\label{Mao10 Conjecture a}
N(0,10,5n) + N(1,10,5n) & > N(4,10,5n) + N(5,10,5n) \text{ for } n \geq 0,\\
\label{Mao10 Conjecture b}
N(1,10,5n) + N(2,10,5n) & \geq N(3,10,5n) + N(4,10,5n) \text{ for } n \geq 1.
\end{align}
\end{conjecture}

The second, is for the $M_2$ rank on partitions without repeated odd parts. 
\begin{conjecture}
Computational evidence suggests that
\begin{align}
\label{Mao610 Conjecture b}
N_2(0,10,5n) + N_2(1,10,5n) & > N_2(4,10,5n) + N_2(5,10,5n) \text{ for } n \geq 0,\\
\label{Mao610 Conjecture c}
N_2(0,10,5n + 4) + N_2(1,10,5n + 4) & > N_2(4,10,5n + 4) + N_2(5,10,5n + 4) \text{ for } n \geq 0,\\
\label{Mao610 Conjecture d}
N_2(1,10,5n) + N_2(2,10,5n) & > N_2(3,10,5n) + N_2(4,10,5n) \text{ for } n \geq 1,\\
\label{Mao610 Conjecture e}
N_2(1,10,5n + 2) + N_2(2,10,5n + 2) & > N_2(3,10,5n + 2) + N_2(4,10,5n + 2) \text{ for } n \geq 1,\\
\label{Mao610 Conjecture a}
N_2(0,6,3n + 2) + N_2(1,6,3n+2) & > N_2(2,6,3n + 2) + N_2(3,6,3n + 2) \text{ for } n \geq 0.
\end{align}
\end{conjecture}

In this paper we prove the following theorem using elementary techniques.
\begin{theorem}\label{main}
Mao's conjectures \eqref{Mao10 Conjecture a}, \eqref{Mao10 Conjecture b}, \eqref{Mao610 Conjecture b}, and \eqref{Mao610 Conjecture c} are true. In fact, in \eqref{Mao10 Conjecture b}, the strict inequality holds.
\end{theorem}
We note that our method did not suffice to prove the remaining three conjectures, which are still open.

The rest of the paper is organized as follows. In Section \ref{prelim}, we gather some definitions, notation, and lemmas that will be used later.  In Section \ref{proof}, we prove Theorem \ref{main}.  

\section{Preliminaries}\label{prelim}

We use the following standard $q$-series notation.  For $n\in\mathbb{N}$, $a\in\mathbb{C}$, define
\begin{align*}
(a;q)_n &:= \prod_{i=0}^{n-1}(1-aq^{i}) \\
(a;q)_\infty &:= \prod_{i=0}^\infty(1-aq^{i}),
\end{align*}
and also define $(a;q)_0=1$.  As shorthand, write 
\begin{align*}
(a_1, \ldots, a_k;q)_n &:= (a_1;q)_n \cdots (a_k;q)_n\\
(a_1, \ldots, a_k;q)_\infty &:= (a_1;q)_\infty \cdots (a_k;q)_\infty.
\end{align*}

Furthermore, we will make use of the following notation of Mao.\footnote{We note that our definition of $L_{a,b}$ differs from Mao's in that the roles of $a$ and $b$ are reversed.}  For positive integers $a<b$, define
\begin{align*}
J_b &:= (q^b; q^b)_{\infty} \\
J_{a,b} &:= (q^a, q^{b-a}, q^b; q^b)_{\infty},\\
L_{a,b} &:=\frac{J_b^2}{J_{a,b}}. 
\end{align*}
\begin{lemma}[Mao \cite{Mao1}]\label{L lemma}
Given positive integers $a<b$, the $q$-series coefficients of $L_{a,b}$ are all nonnegative.
\end{lemma}

Mao proved rank difference formulas that we will use in our proof of Theorem \ref{main}.  First, for unrestricted partitions, Mao proved the following theorem.
\begin{theorem}[Mao \cite{Mao1}]\label{MaoThm1}
We have that
\begin{align*}
&\sum_{n=0}^{\infty} \big(N(0,10,n) + N(1,10,n) - N(4,10,n) - N(5,10,n)\big)q^n \\
&= \left(\frac{J_{25} J_{50}^5 J_{20,50}^2 }{ J_{10,50}^4J_{15,50}^3} + \frac{1}{J_{25}} \sum_{n=-\infty}^{\infty} \frac{(-1)^n q^{75n(n + 1)/2 + 5}}{1 + q^{25n + 5}}\right)  + q\left(\frac{J_{25} J_{50}^5}{J_{5,50} J_{10,50}^2 J_{15,50}^2}\right) \\
&+ q^2\left(\frac{J_{25} J_{50}^5}{J_{5,50}^2 J_{15,50} J_{20,50}^2}\right) + q^3\left(\frac{J_{25} J_{50}^5J_{10,50}^2 }{J_{5,50}^3 J_{20,50}^4} - \frac{1}{J_{25}} \sum_{n=-\infty}^{\infty} \frac{(-1)^n q^{75n(n + 1)/2 + 5}}{1 + q^{25n + 10}}\right) \\
&+ q^4\left(\frac{2 J_{50}^6}{J_{25} J_{5,50} J_{10,50} J_{15,50} J_{20,50}}\right), \\
\end{align*}
and
\begin{align*}
&\sum_{n=0}^{\infty} \big(N(1,10,n) + N(2,10,n) - N(3,10,n) - N(4,10,n)\big)q^n \\
&= \left(\frac{2q^5 J_{50}^6}{J_{25} J_{10,50}^2 J_{15,50}^2} - \frac{1}{J_{25}} \sum_{n=-\infty}^{\infty} \frac{(-1)^n q^{75n(n + 1)/2 + 5}}{1 + q^{25n + 5}}\right) \\
& + q\left(\frac{2q^5 J_{50}^6}{J_{25} J_{5,50} J_{15,50} J_{20,50}^2}\right) + q^2\left(\frac{J_{25}J_{50}^5 J_{20,50} }{J_{10,50}^3 J_{15,50}^3}\right) + q^3\left(\frac{J_{25} J_{50}^5}{J_{5,50} J_{10,50} J_{15,50}^2 J_{20,50} }\right) \\
& + q^4\left(\frac{J_{25} J_{50}^5J_{20,50}^2 J_{25,50}  }{2q^5 J_{10,50}^4 J_{15,50}^4} - \frac{1}{q^5 J_{25}} \sum_{n=-\infty}^{\infty} \frac{(-1)^n q^{(75n^2 + 25n)/2}}{1 + q^{25n}}\right). \\
\end{align*}
\end{theorem}

We will also make use of the following theorem for $M_2$ rank of partitions without repeated odd parts.

\begin{theorem}[Mao \cite{Mao2}] \label{MaoThm2}
We have that
\begin{align*}
& \sum_{n=0}^{\infty} \big(N_2(0,10,n) + N_2(1,10,n) - N_2(4,10,n) - N_2(5,10,n)\big)q^n \\
&= \left( \frac{2q^5 J_{100}^{15} J_{10,100} J_{50,100} }{J_{5,100}^3 J_{15,100}^2 J_{20,100}^2 J_{25,100}^3 J_{30,100} J_{35,100}^2   J_{45,100}^3} + \frac{1}{J_{25,100}} \sum_{n=-\infty}^{\infty} \frac{(-1)^n q^{50n^2 + 25n}}{1 + q^{50n + 10}}\right) \\
& + q \left(\frac{J_{100}^{15} J_{20,100} J_{30,100}^2 J_{50,100}}{ J_{5,100}^2 J_{10,100}^2 J_{15,100}^4 J_{25,100}J_{35,100}^4 J_{40,100}^3J_{45,100}^2  } \right) \\
&+ q^2 \left(\frac{J_{100}^{15} J_{50,100} }{ J_{5,100}^3 J_{15,100}^3 J_{20,100} J_{25,100}J_{35,100}^3 J_{40,100} J_{45,100}^3}\right) \\
&+ q^3\left(\frac{J_{100}^{15}J_{10,100}^2J_{40,100} J_{50,100}  }{ J_{5,100}^4J_{15,100}^2  J_{20,100}^3 J_{25,100} J_{30,100}^2 J_{35,100}^2J_{45,100}^4}\right) \\
& + q^4\left(\frac{2J_{100}^{15}J_{30,100} J_{50,100} }{J_{5,100}^2 J_{10,100} J_{15,100}^3 J_{25,100}^3 J_{35,100}^3J_{40,100}^2 J_{45,100}^2 } + \frac{1}{J_{25,100}} \sum_{n=-\infty}^{\infty} \frac{(-1)^n q^{50n^2 + 75n +20}}{1 + q^{50n + 30}}\right).
\end{align*}
\end{theorem}

In addition, we require the following two facts about $q$-series which follow directly from the definitions.  For integers $a,b,c$ we have
\begin{align}
\label{lemma1} 
(q^a;q^b)_\infty(-q^a;q^b)_\infty & = (q^{2a};q^{2b})_\infty, \\ 
\label{lemma2}
(cq^a;q^{2b})_\infty(cq^{a+b};q^{2b})_\infty & = (cq^{a};q^{b})_\infty. 
\end{align}

Finally, we recall the Jacobi Triple Product formula, which can be found in \cite{Andrews},
\begin{equation}\label{JTP}
\sum_{n \in \mathbb{Z}} z^n q^{n^2} = (-zq,-q/z,q^2;q^2)_{\infty}.
\end{equation}

\section{Proof of Theorem \ref{main}}\label{proof}

\subsection{Proof of (\ref{Mao10 Conjecture a})}\label{first proof}

In order to prove (\ref{Mao10 Conjecture a}), we need to show that the series
\[
\sum_{n=0}^\infty (N(0,10,5n) + N(1,10,5n) - N(4,10,5n) - N(5,10,5n))q^n
\]
has strictly positive coefficients.  Using the first part of Theorem \ref{MaoThm1}, we see that 
\[
\sum_{n=0}^{\infty} \big(N(0,10,n) + N(1,10,n) - N(4,10,n) - N(5,10,n)\big)q^n= S_0 + qS_1 + q^2S_2 + q^3S_3 + q^4S_4, 
\]
where each $S_i$ is a series in $q^5$.  Thus we can obtain our desired generating function by letting $q\mapsto q^\frac{1}{5}$ in $S_0$.  We obtain that
\begin{multline}\label{gen1}
\sum_{n=0}^\infty (N(0,10,5n) + N(1,10,5n) - N(4,10,5n) - N(5,10,5n))q^n \\
= \frac{J_5 J_{10}^5 J_{4,10}^2 }{J_{2,10}^4 J_{3,10}^3  } + \frac{1}{J_5} \sum_{n=-\infty}^{\infty} \frac{(-1)^n q^{(15n^2 + 15n)/2 + 1}}{1 + q^{5n + 1}} = \frac{1}{J_5}\left(\frac{J_5^2 J_{10}^5 J_{4,10}^2 }{J_{2,10}^4 J_{3,10}^3 } + \sum_{n=-\infty}^{\infty} \frac{(-1)^n q^{(15n^2 + 15n)/2 + 1}}{1 + q^{5n + 1}}\right).
\end{multline}

We will now show that \eqref{gen1} has strictly positive $q$-series coefficients for $n \geq 0$. Since $\frac{1}{J_5}$ has all nonnegative coefficients and a constant term of $1$, it suffices to show that 
\[
\frac{J_5^2 J_{10}^5 J_{4,10}^2 }{J_{2,10}^4 J_{3,10}^3 } + \sum_{n=-\infty}^{\infty} \frac{(-1)^n q^{(15n^2 + 15n)/2 + 1}}{1 + q^{5n + 1}}
\]
 has all positive coefficients.  First, we split the sum into nonnegative and negative indices, and reindex to see that
\begin{multline*}
 \sum_{n=-\infty}^{\infty} \frac{(-1)^n q^{(15n^2 + 15n)/2 + 1}}{1 + q^{5n + 1}}  = \sum_{n=0}^{\infty} \frac{(-1)^n q^{(15n^2 + 15n)/2 + 1}}{1 + q^{5n + 1}} + \sum_{n=1}^{\infty} \frac{(-1)^n q^{(15n^2 - 5n)/2}}{1 + q^{5n - 1}} \\
= \sum_{n=0}^{\infty} \frac{(-1)^n q^{(15n^2 + 15n)/2 + 1}(1 - q^{5n + 1})}{1 - q^{10n + 2}} + \sum_{n=1}^{\infty} \frac{(-1)^n q^{(15n^2 - 5n)/2}(1 - q^{5n - 1})}{1 - q^{10n - 2}}. 
\end{multline*}
Now, we split according to the summation index $n$ modulo $2$, to obtain
\begin{multline*}
\sum_{n=-\infty}^{\infty} \frac{(-1)^n q^{(15n^2 + 15n)/2 + 1}}{1 + q^{5n + 1}} \\
= \sum_{n=0}^{\infty} \frac{q^{(15(2n)^2 + 15(2n))/2 + 1}(1 - q^{5(2n) + 1})}{1 - q^{10(2n) + 2}} - \sum_{n=0}^{\infty} \frac{q^{(15(2n+1)^2 + 15(2n+1))/2 + 1}(1 - q^{5(2n+1) + 1})}{1 - q^{10(2n+1) + 2}} \\
+ \sum_{n=1}^{\infty} \frac{q^{(15(2n)^2 - 5(2n))/2}(1 - q^{5(2n) - 1})}{1 - q^{10(2n) - 2}} - \sum_{n=1}^{\infty} \frac{q^{(15(2n-1)^2 - 5(2n-1))/2}(1 - q^{5(2n-1) - 1})}{1 - q^{10(2n-1) - 2}}. 
\end{multline*}

Gathering the positive summands together, we see that
\[
 \sum_{n=-\infty}^{\infty} \frac{(-1)^n q^{(15n^2 + 15n)/2 + 1}}{1 + q^{5n + 1}} = S - T_1 - T_2 - T_3 - T_4,
\]
where
\[
S := \sum_{n=0}^{\infty}\frac{q^{30n^2 +15n + 1}}{1 - q^{20n + 2}} + \sum_{n=0}^{\infty}\frac{q^{30n^2+55n+22}}{1 - q^{20n+12}} + \sum_{n=1}^{\infty}\frac{q^{30n^2-5n}}{1 - q^{20n-2}}+ \sum_{n=1}^{\infty}\frac{q^{30n^2 -25n+ 4 }}{1 - q^{20n-12}},
\]
and
\begin{align*}
T_1 = \sum_{n=0}^\infty a_1(n)q^n &:= \sum_{n=0}^{\infty}\frac{q^{30n^2+25n + 2}}{1 - q^{20n + 2}}, \\
T_2 = \sum_{n=0}^\infty a_2(n)q^n &:= \sum_{n=0}^{\infty}\frac{q^{30n^2+45n+16}}{1 - q^{20n+12}}, \\
T_3 = \sum_{n=0}^\infty a_3(n)q^n &:= \sum_{n=1}^{\infty}\frac{q^{30n^2 + 5n - 1}}{1 - q^{20n-2}}, \\
T_4 = \sum_{n=0}^\infty a_4(n)q^n &:= \sum_{n=1}^{\infty}\frac{q^{30n^2 -35n+10}}{1 - q^{20n-12}}. \\
\end{align*}
We see that $S$, $T_1,\ldots,T_4$ all have nonnegative coefficients.  Thus to prove (\ref{Mao10 Conjecture a}), it suffices to show that
\[
\frac{J_5^2 J_{10}^5 J_{4,10}^2 }{J_{3,10}^3 J_{2,10}^4} - T_1 - T_2 - T_3 - T_4
\]
has positive coefficients.  Let $T_1+T_2+T_3+T_4 = \sum_{n=1}^\infty a(n)q^n$, and let 
\[
\frac{J_5^2 J_{10}^5 J_{4,10}^2 }{J_{3,10}^3 J_{2,10}^4} = 1 +\sum_{n=1}^\infty b(n)q^n.
\] 
We will show that $b(n)>a(n)$ for all $n\geq 1$.

Expanding the denominator of $T_1$ as a geometric series, we see that
\[
T_1 = \sum_{n=0}^\infty \sum_{k=0}^\infty q^{30n^2 + (20k+25)n + (2k+2)}.
\]
Thus for a given $N\geq 0$, we see that $a_1(N)$ counts the number of nonnegative integer pairs $(n,k)$ such that 
\begin{equation}\label{T1eqn}
N=30n^2 + (20k+25)n + (2k+2).
\end{equation}
Clearly for each choice of $n\geq 0$ there is at most one $k\geq 0$ such that $(n,k)$ is a solution to \eqref{T1eqn}.  Also, since $(20k+25)n + (2k+2)$ is positive for all $n,k\geq 0$, if $n\geq \sqrt{\frac{N}{30}}$, then no solutions exist.  Thus, we have that $a_1(N)\leq \lfloor \sqrt{\frac{N}{30}} \rfloor +1$ for all $N\geq 0$.  

Similarly, 
\[
T_2 =  \sum_{n=0}^\infty \sum_{k=0}^\infty q^{30n^2 + (20k+45)n + (12k+16)},
\]
and so $a_2(N)\leq \lfloor \sqrt{\frac{N}{30}}\rfloor +1$ for all $N\geq 0$ as well.  For $T_3$, we have
\[
T_3 =  \sum_{n=1}^\infty \sum_{k=0}^\infty q^{30n^2 + (20k+5)n - (2k+1)}.
\]
Since the sum starts at $n=1$ we have one fewer term.  Also, we see that $(20k+5)n - (2k+1)$ is positive for all $n\geq1,k\geq0$.  Thus we get a bound of $a_3(N)\leq \lfloor \sqrt{\frac{N}{30}}\rfloor $ for all $N\geq 1$.  The $T_4$ case is a little different.  Here,
\[
T_4 = \sum_{n=1}^\infty \sum_{k=0}^\infty q^{30n^2 + (20k-35)n - (12k-10)}.
\]
We observe that $30n^2 + (20k-35)n - (12k-10) \geq 5(3n-2)(2n-1) \geq 5(2n-1)^2$ for all $n\geq 1$.  Thus if $2n-1 \geq \sqrt{\frac{N}{5}}$, then no solutions exist to $N=30n^2 + (20k-35)n - (12k-10)$.  We then get a bound of $a_4(N)\leq \lfloor \sqrt{\frac{N}{20}}\rfloor +1$ for all $N\geq 1$.

Together, noting that none of the $T_i$ have a constant term, we see that for any $n\geq 1$,
\begin{equation}\label{a_bound}
a(n) \leq 3\left\lfloor \sqrt{\frac{n}{30}} \right\rfloor + \left\lfloor \sqrt{\frac{n}{20}} \right\rfloor + 3.
\end{equation}

By \eqref{lemma1} and \eqref{lemma2}, we see that 
\begin{align*}
\frac{J_5^2 J_{10}^5 J_{4,10}^2 }{J_{2,10}^4 J_{3,10}^3 } &= \frac{(q^5; q^5)_{\infty}^2 (q^4, q^6; q^{10})_{\infty}^2}{(q^3, q^7; q^{10})_{\infty}^3 (q^2, q^8; q^{10})_{\infty}^4} \\
&= \frac{(q^5; q^5)_{\infty}^2 (-q^2, -q^3, -q^7, -q^8; q^{10})_{\infty}^2}{(q^3, q^7; q^{10})_{\infty} (q^2, q^8; q^{10})_{\infty}^2} \\
&= \frac{(q^5; q^5)_{\infty}^2 (-q^2, -q^3; q^5)_{\infty}^2}{(q^3, q^7; q^{10})_{\infty} (q^2, q^8; q^{10})_{\infty}^2}.
\end{align*}
 
\noindent Applying \eqref{JTP} with $z = q^{1/2}$ and $q = q^{5/2}$, we obtain
\begin{multline*}
\frac{J_5^2 J_{4,10}^2 J_{10}^5}{J_{3,10}^3 J_{2,10}^4} = \frac{1}{(q^3, q^7; q^{10})_{\infty} (q^2, q^8; q^{10})_{\infty}^2} \left[\sum_{n=-\infty}^{\infty} q^{n(5n + 1)/2}\right]^2 \\
= \frac{1}{(1 - q^2)(1 - q^3)} \cdot \frac{1}{(q^2, q^7, q^{12}, q^{13}; q^{10})_{\infty} (q^8; q^{10})_{\infty}^2}\left[\sum_{n=-\infty}^{\infty} q^{n(5n + 1)/2}\right]^2,
\end{multline*}
where we observe that all series involved in this product have nonnegative coefficients.  Let 
\[
\frac{1}{(1 - q^2)(1 - q^3)} = \sum_{i,j\geq 0} q^{2i+3j}=\sum_{n=0}^{\infty} c(n)q^n.
\]

We note that $c(0)=1$ and for $n\geq 1$, $c(n)$ is equal to the number of nonnegative integer solutions $(i,j)$ of the equation $2i + 3j = n$.  For a fixed $n\geq 1$, and $j\geq 0$, we see that there is at most one $i\geq 0$ for which $(i,j)$ is a solution, and such an $i$ exists if and only if $0\leq j \leq n/3$ and $j\equiv n \pmod 2$.  Considering each possible residue of $n$ modulo $6$, we see that in all cases, $c(n)\geq \lfloor \frac{n}{6} \rfloor $.  Thus, we have that for all $n\geq 1$,
\begin{equation}\label{b_bound}
b(n) \geq \left\lfloor \frac{n}{6} \right\rfloor.
\end{equation} 

It suffices then to show that $\frac{n}{6} > 3\sqrt{\frac{n}{30}}+\sqrt{\frac{n}{20}}  + 4$ for sufficiently large $n$, and to check that $b(n)>a(n)$ for all remaining cases.  We have that $\frac{n}{6} > 3\sqrt{\frac{n}{30}}+\sqrt{\frac{n}{20}}  + 4$ if and only if $\frac{1}{6}n - (\frac{\sqrt{30}+\sqrt{5}}{10})\sqrt{n} -4 \geq 0$, which occurs for $n\geq 60$.  Moreover, we also see that $b(n)>a(n)$ for $1\leq n \leq 59$, by a quick Maple calculation, which completes the proof of \eqref{Mao10 Conjecture a}.

For the remaining conjectures we use a similar technique, so give a somewhat abbreviated discussion of the proofs.

\subsection{Proof of (\ref{Mao10 Conjecture b})}\label{second proof}
In order to prove (\ref{Mao10 Conjecture b}), we need to show that 
\[
\sum_{n=1}^\infty (N(1,10,5n) + N(2,10,5n) - N(3,10,5n) - N(4,10,5n))q^n
\]
has nonnegative coefficients.  Using the second part of Theorem \ref{MaoThm1}, we obtain that 
\begin{multline}\label{gen2}
\sum_{n=1}^\infty (N(1,10,5n) + N(2,10,5n) - N(3,10,5n) - N(4,10,5n))q^n \\
= \frac{1}{J_5} \left(\frac{2qJ_{10}^6}{J_{2,10}^2 J_{3,10}^2} - \sum_{n=-\infty}^{\infty} \frac{(-1)^n q^{(15n^2 + 15n)/2 + 1}}{1 + q^{5n + 1}}\right).
\end{multline}
Since $\frac{1}{J_5}$ has all nonnegative coefficients and a constant term of $1$, it suffices to show that 
\[
\frac{2qJ_{10}^6}{J_{2,10}^2 J_{3,10}^2} - \sum_{n=-\infty}^{\infty} \frac{(-1)^n q^{(15n^2 + 15n)/2 + 1}}{1 + q^{5n + 1}}
\]
has all nonnegative coefficients.  We observe the sum in this case is the same as the sum in the proof of \eqref{Mao10 Conjecture a}.  However in this setting we are subtracting, rather than adding the sum.  Thus by our dissection in the last subsection, if suffices to prove that 
\[
\frac{2qJ_{10}^6}{J_{2,10}^2 J_{3,10}^2} - T_1' - T_2' - T_3' - T_4'
\]
has positive coefficients, where
\begin{align*}
T_1' = \sum_{n=0}^\infty a_1'(n)q^n &:= \sum_{n=0}^{\infty}\frac{q^{30n^2 +15n + 1}}{1 - q^{20n + 2}}= \sum_{n=0}^{\infty}\sum_{k=0}^{\infty} q^{30n^2 + 15n+1 +k(20n+2)}, \\
T_2' = \sum_{n=0}^\infty a_2'(n)q^n &:= \sum_{n=0}^{\infty}\frac{q^{30n^2+55n+22}}{1 - q^{20n+12}}= \sum_{n=0}^{\infty}\sum_{k=0}^{\infty} q^{30n^2 + 55n+22 +k(20n+12)}, \\
T_3' = \sum_{n=0}^\infty a_3'(n)q^n &:= \sum_{n=1}^{\infty}\frac{q^{30n^2-5n}}{1 - q^{20n-2}}= \sum_{n=1}^{\infty}\sum_{k=0}^{\infty} q^{30n^2 -5n +k(20n-2)}, \\
T_4' = \sum_{n=0}^\infty a_4'(n)q^n &:= \sum_{n=1}^{\infty}\frac{q^{30n^2 -25n+ 4 }}{1 - q^{20n-12}}= \sum_{n=1}^{\infty}\sum_{k=0}^{\infty} q^{30n^2 -25n+ 4  +k(20n-12)}. \\
\end{align*}
Let $T_1'+T_2'+T_3'+T_4' = \sum_{n=1}^\infty a'(n)q^n$, and let 
\[
\frac{2qJ_{10}^6}{J_{2,10}^2 J_{3,10}^2} = \sum_{n=1}^\infty b'(n)q^n.
\] 
We will show that $b'(n)> a'(n)$ for all $n\geq 1$.

Arguing as in Section \ref{first proof}, we see that for any $N\geq 1$, $a_1'(N),a_2'(N) \leq \lfloor \sqrt{\frac{N}{30}} \rfloor +1$.  Also, we note that since $30n^2 -5n+k(20n-2) > (5n)^2$ for all $n\geq 1$, we have that $a_3'(N) \leq  \lfloor \sqrt{\frac{N}{25}} \rfloor$.  Similarly, since $30n^2 -25n+ 4  > 30(n-1)^2$ for all $n\geq 1$, we have that $a_4'(N) \leq \lfloor \sqrt{\frac{N}{30}} \rfloor + 1$.  Together, noting that none of the $T_i'$ have a constant term, we see that for any $n\geq 1$,
\begin{equation}\label{a'_bound}
a'(n) \leq 3\left\lfloor \sqrt{\frac{n}{30}} \right\rfloor + \left\lfloor \sqrt{\frac{n}{25}} \right\rfloor + 3.
\end{equation}

We now examine the product. By \eqref{lemma1}, \eqref{lemma2}, and\eqref{JTP}, we see that
\begin{align*}
\frac{2q J_{10}^6}{J_{2,10}^2 J_{3,10}^2} &= 2qL_{3,10} \left(\frac{(-q^2, -q^8; q^{10})_{\infty} (q^{10}; q^{10})_{\infty}}{(q^4, q^{16}; q^{20})_{\infty} (q^2, q^3, q^7, q^8; q^{10})_{\infty}}\right) \\
&= \frac{2q}{(1 - q^2)(1 - q^3)} \left(\frac{L_{3,10}}{(q^4, q^{16}; q^{20})_{\infty} (q^7, q^8, q^{12}, q^{13}; q^{10})_{\infty}} \cdot \sum_{n=-\infty}^{\infty} q^{5n^2 + 2n}\right).
\end{align*}
Arguing as before, we find that the coefficient of $q^n$ in $\frac{2q}{(1 - q^2)(1 - q^3)}$ is at least $2\left\lfloor\frac{n - 1}{6}\right\rfloor$ for $n \geq 1$.  We see that $L_{3,10}$ has a constant term of $1$, and by Lemma \ref{L lemma}, $L_{3,10}$ has all nonnegative coefficients.  Thus we have that for all $n\geq 1$,
\[
b'(n)\geq 2\left\lfloor\frac{n - 1}{6}\right\rfloor.
\]

Since $2\left\lfloor{\frac{n - 1}{6}}\right\rfloor > 3\left\lfloor \sqrt{\frac{n}{30}} \right\rfloor + \left\lfloor \sqrt{\frac{n}{25}} \right\rfloor + 3$ for $n \geq 24$, it thus suffices to check that $b'(n)> a'(n)$ for $1 \leq n \leq 23$. A quick computation in Maple verifies that this is true.

\subsection{Proof of (\ref{Mao610 Conjecture b})}\label{third proof}
In order to prove (\ref{Mao610 Conjecture b}), we need to show that 
\[
\sum_{n=0}^\infty (N_2(0,10,5n) + N_2(1,10,5n) - N_2(4,10,5n) - N_2(5,10,5n))q^n
\]
has positive coefficients.  Using the Theorem \ref{MaoThm2}, we obtain that 
\begin{multline}\label{gen3}
\sum_{n=0}^\infty (N_2(0,10,5n) + N_2(1,10,5n) - N_2(4,10,5n) - N_2(5,10,5n))q^n \\
= \frac{1}{J_{5,20}} \left(\frac{2qJ_{20}^{15}J_{2,20}J_{10,20}}{J_{1,20}^3 J_{3,20}^2J_{4,20}^2J_{5,20}^2J_{6,20}J_{7,20}^2J_{9,20}^3} + \sum_{n=-\infty}^{\infty} \frac{(-1)^n q^{10n^2 + 5n}}{1 + q^{10n + 2}}\right).
\end{multline}
Since $\frac{1}{J_{5,20}}$ has all nonnegative coefficients and a constant term of $1$, it suffices to show that 
\[
\frac{2qJ_{20}^{15}J_{2,20}J_{10,20}}{J_{1,20}^3 J_{3,20}^2J_{4,20}^2J_{5,20}^2J_{6,20}J_{7,20}^2J_{9,20}^3} + \sum_{n=-\infty}^{\infty} \frac{(-1)^n q^{10n^2 + 5n}}{1 + q^{10n + 2}}
\]
has all positive coefficients.  Splitting up the sum as we do in the proof of \eqref{Mao10 Conjecture a}, we obtain that
\[
 \sum_{n=-\infty}^{\infty} \frac{(-1)^n q^{10n^2 + 5n}}{1 + q^{10n + 2}} = S'' - T_1'' - T_2'' - T_3'' - T_4'',
\]
where
\[
S'':= \sum_{n=0}^{\infty}\frac{q^{40n^2 +10n}}{1 - q^{40n + 4}} + \sum_{n=0}^{\infty}\frac{q^{40n^2+70n+27}}{1 - q^{40n+24}} + \sum_{n=1}^{\infty}\frac{q^{40n^2+10n-2}}{1 - q^{40n-4}}+ \sum_{n=1}^{\infty}\frac{q^{40n^2 -10n-9}}{1 - q^{40n-24}},
\]
and
\begin{align*}
T_1'' = \sum_{n=0}^\infty a_1''(n)q^n &:= \sum_{n=0}^{\infty}\frac{q^{40n^2+30n + 2}}{1 - q^{40n + 4}} = \sum_{n=0}^{\infty}\sum_{k=0}^{\infty} q^{40n^2 + 30n+2 +k(40n+4)}, \\
T_2'' = \sum_{n=0}^\infty a_2''(n)q^n &:= \sum_{n=0}^{\infty}\frac{q^{40n^2+50n+15}}{1 - q^{40n+24}}= \sum_{n=0}^{\infty}\sum_{k=0}^{\infty} q^{40n^2 + 50n+15 +k(40n+24)}, \\
T_3'' = \sum_{n=0}^\infty a_3''(n)q^n &:= \sum_{n=1}^{\infty}\frac{q^{40n^2 + 30n - 4}}{1 - q^{40n-4}}= \sum_{n=1}^{\infty}\sum_{k=0}^{\infty} q^{40n^2 + 30n-4 +k(40n-4)}, \\
T_4'' = \sum_{n=0}^\infty a_4''(n)q^n &:= \sum_{n=1}^{\infty}\frac{q^{40n^2 -30n+3}}{1 - q^{40n-24}}= \sum_{n=1}^{\infty}\sum_{k=0}^{\infty} q^{40n^2 -30n+3 +k(40n-24)}. \\
\end{align*}
We see that $S''$, $T_1'',\ldots,T_4''$ all have nonnegative coefficients.  Thus to prove (\ref{Mao610 Conjecture b}), it suffices to show that
\[
\frac{2qJ_{20}^{15}J_{2,20}J_{10,20}}{J_{1,20}^3 J_{3,20}^2J_{4,20}^2J_{5,20}^2J_{6,20}J_{7,20}^2J_{9,20}^3} - T_1'' - T_2'' - T_3'' - T_4''
\]
has positive coefficients. Let $T_1''+T_2''+T_3''+T_4'' = \sum_{n=1}^\infty a''(n)q^n$, and let 
\[
\frac{2qJ_{20}^{15}J_{2,20}J_{10,20}}{J_{1,20}^3 J_{3,20}^2J_{4,20}^2J_{5,20}^2J_{6,20}J_{7,20}^2J_{9,20}^3} = \sum_{n=1}^\infty b''(n)q^n.
\] 
We will show that $b''(n) > a''(n)$ for all $n\geq 1$.

Again arguing as in Section \ref{first proof}, we see that for any $N\geq 1$, $a_1''(N),a_2''(N) \leq \lfloor \sqrt{\frac{N}{40}} \rfloor +1$, and $a_3''(N) \leq  \lfloor \sqrt{\frac{N}{40}} \rfloor$.  Also, since $40n^2 -30n+3 > 6(2n-1)^2$ for all $n\geq 1$, we have that $a_4''(N) \leq \lfloor \sqrt{\frac{N}{24}} \rfloor + 1$.  Together, noting that none of the $T_i''$ have a constant term, we see that for any $n\geq 1$,
\begin{equation}\label{a''_bound}
a''(n) \leq 3\left\lfloor \sqrt{\frac{n}{40}} \right\rfloor + \left\lfloor \sqrt{\frac{n}{24}} \right\rfloor + 3.
\end{equation}

We now examine the product. By \eqref{lemma1}, we find that
\begin{multline*}
\frac{2qJ_{2,20} J_{10,20} J_{20}^{15}}{J_{6,20}J_{3,20}^2 J_{4,20}^2 J_{5,20}^2 J_{7,20}^2 J_{1,20}^3 J_{9,20}^3} \\
= \frac{2q}{(1-q)^2}L_{9,20}^2 \left(\frac{(-q, -q^9, -q^{11}, -q^{19}; q^{20})_{\infty} (-q^5, -q^{15}; q^{20})_{\infty}^2}{(q^6, q^{14}; q^{20})_{\infty} (q^3, q^4, q^7, q^{13}, q^{16}, q^{17}, q^{19},q^{21}; q^{20})_{\infty}^2 (q^{19}; q^{20})_{\infty}^3}\right). \\
\end{multline*}
Expanding gives that $\frac{2q}{(1-q)^2} = \sum_{n=1}^\infty 2nq^n$.  We know by Lemma \ref{L lemma} that $L_{9,20}$ has nonnegative coefficients and a constant term of $1$, and we can observe that this is true for the rest of the product as well.  Thus, we have for all $n\geq 1$, 
\[
b''(n)\geq 2n.
\]
Since $2n > 3\left\lfloor \sqrt{\frac{n}{40}} \right\rfloor + \left\lfloor \sqrt{\frac{n}{24}} \right\rfloor + 3$ for $n \geq 2$, it thus suffices to observe that $2=b''(1) > a''(1)=0$.

\subsection{Proof of (\ref{Mao610 Conjecture c})}\label{fourth proof}
In order to prove (\ref{Mao610 Conjecture c}), we need to show that 
\[
\sum_{n=0}^\infty (N_2(0,10,5n+4) + N_2(1,10,5n+4) - N_2(4,10,5n+4) - N_2(5,10,5n+4))q^n
\]
has positive coefficients.  Using the Theorem \ref{MaoThm2}, we obtain that 
\begin{multline}\label{gen4}
\sum_{n=0}^\infty (N_2(0,10,5n+4) + N_2(1,10,5n+4) - N_2(4,10,5n+4) - N_2(5,10,5n+4))q^n \\
= \frac{1}{J_{5,20}} \left(\frac{2J_{20}^{15}J_{6,20}J_{10,20}}{J_{1,20}^2 J_{2,20}J_{3,20}^3J_{5,20}^2J_{7,20}^3J_{8,20}^2J_{9,20}^2} + \sum_{n=-\infty}^{\infty} \frac{(-1)^n q^{10n^2 + 15n+4}}{1 + q^{10n + 6}}\right).
\end{multline}
Since $\frac{1}{J_{5,20}}$ has all nonnegative coefficients and a constant term of $1$, it suffices to show that 
\[
\frac{2J_{20}^{15}J_{6,20}J_{10,20}}{J_{1,20}^2 J_{2,20}J_{3,20}^3J_{5,20}^2J_{7,20}^3J_{8,20}^2J_{9,20}^2} + \sum_{n=-\infty}^{\infty} \frac{(-1)^n q^{10n^2 + 15n+4}}{1 + q^{10n + 6}}
\]
has all positive coefficients.  Splitting up the sum as we do in the proof of \eqref{Mao10 Conjecture a}, we obtain that
\[
 \sum_{n=-\infty}^{\infty} \frac{(-1)^n q^{10n^2 + 15n +4}}{1 + q^{10n + 6}} = S''' - T_1''' - T_2''' - T_3''' - T_4''',
\]
where
\[
S''':= \sum_{n=0}^{\infty}\frac{q^{40n^2 +30n+4}}{1 - q^{40n + 12}} + \sum_{n=0}^{\infty}\frac{q^{40n^2+ 90n+45 }}{1 - q^{40n+32}} + \sum_{n=1}^{\infty}\frac{q^{40n^2-10n-2}}{1 - q^{40n-12}}+ \sum_{n=1}^{\infty}\frac{q^{40n^2  -30n-3 }}{1 - q^{40n-32}},
\]
and
\begin{align*}
T_1''' = \sum_{n=0}^\infty a_1'''(n)q^n &:= \sum_{n=0}^{\infty}\frac{q^{40n^2+50n + 10}}{1 - q^{40n + 12}} = \sum_{n=0}^{\infty}\sum_{k=0}^{\infty} q^{40n^2+50n + 10 +k(40n+12)}, \\
T_2''' = \sum_{n=0}^\infty a_2'''(n)q^n &:= \sum_{n=0}^{\infty}\frac{q^{40n^2+ 70n+29 }}{1 - q^{40n+32}}= \sum_{n=0}^{\infty}\sum_{k=0}^{\infty} q^{40n^2+ 70n+29  +k(40n+32)}, \\
T_3''' = \sum_{n=0}^\infty a_3'''(n)q^n &:= \sum_{n=1}^{\infty}\frac{q^{40n^2 + 10n - 8}}{1 - q^{40n-12}}= \sum_{n=1}^{\infty}\sum_{k=0}^{\infty} q^{40n^2 + 10n - 8 +k(40n-12)}, \\
T_4''' = \sum_{n=0}^\infty a_4'''(n)q^n &:= \sum_{n=1}^{\infty}\frac{q^{40n^2 -50n+13}}{1 - q^{40n-32}}= \sum_{n=1}^{\infty}\sum_{k=0}^{\infty} q^{40n^2  -50n+13  +k(40n-32)}. \\
\end{align*}
We see that $S'''$, $T_1''',\ldots,T_4'''$ all have nonnegative coefficients.  Thus to prove (\ref{Mao610 Conjecture c}), it suffices to show that
\[
\frac{2J_{20}^{15}J_{6,20}J_{10,20}}{J_{1,20}^2 J_{2,20}J_{3,20}^3J_{5,20}^2J_{7,20}^3J_{8,20}^2J_{9,20}^2} - T_1''' - T_2''' - T_3''' - T_4'''
\]
has positive coefficients. Let $T_1'''+T_2'''+T_3'''+T_4''' = \sum_{n=1}^\infty a'''(n)q^n$, and let 
\[
\frac{2J_{20}^{15}J_{6,20}J_{10,20}}{J_{1,20}^2 J_{2,20}J_{3,20}^3J_{5,20}^2J_{7,20}^3J_{8,20}^2J_{9,20}^2} = 2+ \sum_{n=1}^\infty b'''(n)q^n.
\] 
We will show that $b'''(n) > a'''(n)$ for all $n\geq 1$.

Again arguing as in Section \ref{first proof}, we see that for any $N\geq 1$, $a_1'''(N),a_2'''(N) \leq \lfloor \sqrt{\frac{N}{40}} \rfloor +1$, and $a_3'''(N) \leq  \lfloor \sqrt{\frac{N}{40}} \rfloor$.  Also, since $40n^2 -50n+13 > 40(n-1)^2$ for all $n\geq 1$, we have that $a_4'''(N) \leq \lfloor \sqrt{\frac{N}{40}} \rfloor + 1$.  Together, noting that none of the $T_i'''$ have a constant term, we see that for any $n\geq 1$,
\begin{equation}\label{a'''_bound}
a'''(n) \leq 4\left\lfloor \sqrt{\frac{n}{40}} \right\rfloor + 3.
\end{equation}

By computing the $n = 0$ term of the sum appearing in (\ref{Mao610 Conjecture c}), we find that the constant term is $1$. We may thus instead consider

\begin{equation}
\label{Mao610c Reformulated}
\frac{2J_{6,20} J_{10,20} J_{20}^{15}}{J_{2,20} J_{1,20}^2 J_{5,20}^2 J_{8,20}^2 J_{9,20}^2 J_{3,20}^3 J_{7,20}^3} - \frac{1}{(1 - q)^2}.
\end{equation}

We now examine the product. By \eqref{lemma1}, we find that
\begin{multline*}
\frac{2 J_{6,20} J_{10,20} J_{100}^{15}}{J_{2,20} J_{1,20}^2 J_{5,20}^2 J_{8,20}^2 J_{9,20}^2 J_{3,20}^3 J_{7,20}^3} \\
= \frac{2}{(1-q)^2}L_{9,20}^2 \left(\frac{(-q^{3}, -q^{7}, -q^{13}, -q^{17}; q^{20})_{\infty} (-q^{5}, -q^{15}; q^{20})_{\infty}^2}{(q^2, q^{18}; q^{20})_{\infty} (q^3, q^7, q^8, q^{12}, q^{13}, q^{17}, q^{19}, q^{21}; q^{20})_{\infty}^2}\right). 
\end{multline*}

As in Section \ref{third proof}, we have that expanding gives that $\frac{2}{(1-q)^2} =\sum_{n=0}^\infty 2(n+1)q^n$.  Also, we have already seen that $L_{9,20}$ has nonnegative coefficients and a constant term of $1$, and we can observe that this is true for the rest of the product as well.  Thus, we have for all $n\geq 1$, 
\[
b'''(n)\geq 2(n+1).
\]
Since $2(n+1) > 4\left\lfloor \sqrt{\frac{n}{40}} \right\rfloor + 3$ for $n \geq 1$, this completes the proof of \eqref{Mao610 Conjecture c}. 



\bibliography{REU15}

\begin{thebibliography}{1}

\bibitem{Andrews}
George~E Andrews.
\newblock {\em The theory of partitions}.
\newblock Number~2. Cambridge university press, 1998.

\bibitem{ASD}
A.~O.~L. Atkin and P.~Swinnerton-Dyer.
\newblock Some properties of partitions.
\newblock {\em Proc. London Math. Soc. (3)}, 4:84--106, 1954.

\bibitem{Dyson}
Freeman~J Dyson.
\newblock Some guesses in the theory of partitions.
\newblock {\em Eureka (Cambridge)}, 8(10), 1944.

\bibitem{LO1}
Jeremy Lovejoy and Robert Osburn.
\newblock Rank differences for overpartitions.
\newblock {\em Q. J. Math.}, 59(2):257--273, 2008.

\bibitem{LO3}
Jeremy Lovejoy and Robert Osburn.
\newblock {$M_2$}-rank differences for partitions without repeated odd parts.
\newblock {\em J. Th\'eor. Nombres Bordeaux}, 21(2):313--334, 2009.

\bibitem{LO2}
Jeremy Lovejoy and Robert Osburn.
\newblock {$M_2$}-rank differences for overpartitions.
\newblock {\em Acta Arith.}, 144(2):193--212, 2010.

\bibitem{Mao1}
Renrong Mao.
\newblock Ranks of partitions modulo 10.
\newblock {\em J. Number Theory}, 133(11):3678--3702, 2013.

\bibitem{Mao2}
Renrong Mao.
\newblock The {$M_2$}-rank of partitions without repeated odd parts modulo
  {$6$} and {$10$}.
\newblock {\em Ramanujan J.}, 37(2):391--419, 2015.

\end{thebibliography}
\bibliographystyle{plain}

\end{document}